\documentclass[10pt,twoside]{article} 
\usepackage{amsfonts}
\usepackage{amssymb,amsmath,amsthm,harvard,ecphead}
\citationstyle{agsm}
\newif\ifpdf
\ifx\pdfoutput\undefined
\pdffalse%we are not running pdflatex
\else
\pdfoutput=1%we are running pdflatex
\pdftrue
\fi
\ifpdf
\fi
\makeatletter
\renewcommand\section{\@startsection {section}{1}{\z@}%
                                   {-3.5ex \@plus -1ex \@minus -.2ex}%
                                   {2.3ex \@plus.2ex}%
                                   {\centering\normalfont\normalsize\scshape}}
\makeatother
\begin{document}
% ecp-paper-head.sty
%%%%%%%%%%%%%%%%%%%%%%%%%%%%%%%%%%%%%%%%%%%%%%%%%%%%%%%%%%%%%%%%%%%%%%%%%%%
%
%       ELECTRONIC COMMUNICATION IN PROBABILITY
%       PAPER-SPECIFIC STYLE FILE: CREATES THE AUTHOR NAMES, ADDRESSES, ETC.
%       SHOULD BE SAVE AS `paper?.head.sty', WHERE `?' DENOTES THE PAPER NO.
%
%       AUTHORS: SEARCH FOR `XXX' or `%XXX' AND MAKE THE APPROPRIATE CHANGES
%       Editor/Managing Editor - look for YYY
%
%       Author: Davar Khoshnevisan
%       First Update : February 20, 2000
%       Second Update: June 7, 2000
%       Third update (Martin Barlow) 30 September 2002
%       Updates to use ecp03.sty - 27 Jan 2003, 14 Feb. 2003
%%%%%%%%%%%%%%%%%%%%%%%%%%%%%%%%%%%%%%%%%%%%%%%%%%%%%%%%%%%%%%%%%%%%%%%%%%%
%
%
%               TITLE AND SHORT TITLE
%
%%%%%%%%%%%%%%%%%%%%%%%%%%%%%%%%%%%%%%%%%%%%%%%%%%%%%%%%%%%%
%
\newcommand{\Title}{\uppercase{%
State Tameness: A New Approach for Credit Constrains
}}
\newcommand{\ShortTitle}{%
State Tameness
}
%%%%%%%%%%%%%%%%%%%%%%%%%%%%%%%%%%%%%%%%%%%%%%%%%%%%%%%%%%%%
%
%               AUTHOR and NUMBER OF AUTHORS
%
%%%%%%%%%%%%%%%%%%%%%%%%%%%%%%%%%%%%%%%%%%%%%%%%%%%%%%%%%%%%
%
\newcounter{NumAuthor}
\setcounter{NumAuthor}{%
1 %XXX    NUMBER OF AUTHORS
}
\addtocounter{NumAuthor}{-1}
\newcommand{\AuthorOne}{\uppercase{%
%XXX Fill in AUTHOR 1's name and support if appropriate 
%    - see the 2 examples below
Jaime A. Londo\~no
}}
\newcommand{\AuthorTwo}{\uppercase{%
%XXX    AUTHOR 2
}}
\newcommand{\AuthorThree}{\uppercase{%
%XXX    AUTHOR 3
}}
\newcommand{\AuthorFour}{\uppercase{%
%XXX    AUTHOR 4
}}
\newcommand{\AuthorFive}{\uppercase{%
%XXX    AUTHOR 5
}}
%%%%%%%%%%%%%%%%%%%%%%%%%%%%%%%%%%%%%%%%%%%%%%%%%%%%%%%%%%%%
%
%               ADDRESS
%
%%%%%%%%%%%%%%%%%%%%%%%%%%%%%%%%%%%%%%%%%%%%%%%%%%%%%%%%%%%%
%
\newcommand{\AddressOne}{%
Departamento de Ciencias B\'asicas,\\
Universidad EAFIT
}
\newcommand{\AddressTwo}{%
2st Ave., 23 rd.St., NJ. 84105          %XXX ADDRESS FOR AUTHOR 2
\\United States of America              %XXX (extra line not nec.)
}
\newcommand{\AddressThree}{%
2st Ave., 23 rd.St., NJ. 84105          %XXX ADDRESS FOR AUTHOR 3
\\United States of America              %XXX (extra line not nec.)
}
\newcommand{\AddressFour}{%
2st Ave., 23 rd.St., NJ. 84105          %XXX ADDRESS FOR AUTHOR 4
\\United States of America              %XXX (extra line not nec.)
}
\newcommand{\AddressFive}{%
2st Ave., 23 rd.St., NJ. 84105          %XXX ADDRESS FOR AUTHOR 5
\\United States of America              %XXX (extra line not nec.)
}
%%%%%%%%%%%%%%%%%%%%%%%%%%%%%%%%%%%%%%%%%%%%%%%%%%%%%%%%%%%%
%
%               EMAIL
%
%%%%%%%%%%%%%%%%%%%%%%%%%%%%%%%%%%%%%%%%%%%%%%%%%%%%%%%%%%%%
%
\newcommand{\EmailOne}{%
jalondon@eafit.edu.co
}
\newcommand{\EmailTwo}{%
email.com                               %XXX EMAIL FOR AUTHOR 2
}
\newcommand{\EmailThree}{%
email.com                               %XXX EMAIL FOR AUTHOR 3
}
\newcommand{\EmailFour}{%
email.com                               %XXX EMAIL FOR AUTHOR 4
}
\newcommand{\EmailFive}{%
email.com                               %XXX EMAIL FOR AUTHOR 5
}
%%%%%%%%%%%%%%%%%%%%%%%%%%%%%%%%%%%%%%%%%%%%%%%%%%%%%%%%%%%%
%
%               PAGE INFORMATION
%
%%%%%%%%%%%%%%%%%%%%%%%%%%%%%%%%%%%%%%%%%%%%%%%%%%%%%%%%%%%%
%  Authors - don't worry about these, we will reset these lines.
\setcounter{volume}{9}                  %YYY VOLUME
\setcounter{year}{2004}                     %YYY YEAR
\setcounter{firstpage}{1}               %YYY STARTING PAGE
\setcounter{lastpage}{13}                %YYY LAST PAGE
\setcounter{page}{\value{firstpage}}
%
%%%%%%%%%%%%%%%%%%%%%%%%%%%%%%%%%%%%%%%%%%%%%%%%%%%%%%%%%%%%
%
%               SUBMISSION/ACCEPTANCE DATES
%
%%%%%%%%%%%%%%%%%%%%%%%%%%%%%%%%%%%%%%%%%%%%%%%%%%%%%%%%%%%%
%
\newcommand{\Submitted}{%
9 May 2003
}
\newcommand{\Accepted}{%
3 February 2004
}
%%%%%%%%%%%%%%%%%%%%%%%%%%%%%%%%%%%%%%%%%%%%%%%%%%%%%%%%%%%%
%
%               Subject Classification and Keywords
%
%%%%%%%%%%%%%%%%%%%%%%%%%%%%%%%%%%%%%%%%%%%%%%%%%%%%%%%%%%%%
%
\newcommand{\SubjectClassification}{%
Primary 91B28, 60G40; secondary 60H10.
}
\newcommand{\Keywords}{%
arbitrage, pricing of contingent claims,
  continuous-time financial markets, tameness.
}
%%%%%%%%%%%%%%%%%%%%%%%%%%%%%%%%%%%%%%%%%%%%%%%%%%%%%%%%%%%%
%
%               Abstract
%
%%%%%%%%%%%%%%%%%%%%%%%%%%%%%%%%%%%%%%%%%%%%%%%%%%%%%%%%%%%%
%
\newcommand{\Abstract}{%
 We propose
  a new definition for tameness within the model of security prices as
  It\^o processes that is risk-aware.  We give a new definition for
  arbitrage and characterize it.  We then prove a theorem that can be
  seen as an extension of the second fundamental theorem of asset
  pricing, and a theorem for valuation of contingent claims of the
  American type. The valuation of European contingent claims and
  American contingent claims that we obtain does not require the full
  range of the volatility matrix.  The technique used to prove the
  theorem on valuation of American contingent claims does not depend
  on the Doob-Meyer decomposition of super-martingales; its proof is
  constructive and suggest and alternative way to find approximations of
  stopping times that are close to optimal.
}
%%%%%%%%%%%%%%%%%%%%%%%%%%%%%%%%%%%%%%%%%%%%%%%%%%%%%%%%%%%%%%%
%    ecp03.sty
%
%    ELECTRONIC COMMUNICATION IN PROBABILITY
%    MAIN STYLE FILE
%
%           PLEASE DO NOT EDIT!
%
%    Author: Davar Khoshnevisan
%    February 20, 2000
%    Minor changes to style, (Martin Barlow) 24 January 2003
%    This file (ecp03.sty) for Volume 8 on.
%
%%%%%%%%%%%%%%%%%%%%%%%%%%%%%%%%%%%%%%%%%%%%%%%%%%%%%%%%%%%%%%%
%
\thispagestyle{plain}
\rule{0in}{0pt}

%%%%%%%%%%%%%%%%%%%%%%%%%%%%%%%%%%%%%%%%%%%%%%%%%%%%%%%%%%%%
%
%  Add Title
%
%%%%%%%%%%%%%%%%%%%%%%%%%%%%%%%%%%%%%%%%%%%%%%%%%%%%%%%%%%%%
\vskip 60pt
\noindent
{\Large\bf\uppercase{\Title} }  

%%%%%%%%%%%%%%%%%%%%%%%%%%%%%%%%%%%%%%%%%%%%%%%%%%%%%%%%%%%%
%
%  Add Author, Address, Email
%
%%%%%%%%%%%%%%%%%%%%%%%%%%%%%%%%%%%%%%%%%%%%%%%%%%%%%%%%%%%%
\vskip 12pt\noindent
\ifcase\value{NumAuthor}
        \uppercase{\AuthorOne}\\
        {\em \AddressOne}\\
        email:\hskip 4pt\texttt{\EmailOne}
\or
        \uppercase{\AuthorOne}\\
        {\em \AddressOne}\\
        email:\hskip 4pt\texttt{\EmailOne}\\[2mm]
        \uppercase{\AuthorTwo}\\
        {\em \AddressTwo}\\
        email:\hskip 4pt\texttt{\EmailTwo}
\or
        \uppercase{\AuthorOne}\\
        {\em \AddressOne}\\
        email:\hskip 4pt\texttt{\EmailOne}\\[2mm]
        \uppercase{\AuthorTwo}\\
        {\em \AddressTwo}\\
        email:\hskip 4pt\texttt{\EmailTwo}\\[2mm]
        \uppercase{\AuthorThree}\\
        {\em \AddressThree}\\
        email:\hskip 4pt\texttt{\EmailThree}
\or
        \uppercase{\AuthorOne}\\
        {\em \AddressOne}\\
        email:\hskip 4pt\texttt{\EmailOne}\\[2mm]
        \uppercase{\AuthorTwo}\\
        {\em \AddressTwo}\\
        email:\hskip 4pt\texttt{\EmailTwo}\\[2mm]
        \uppercase{\AuthorThree}\\
        {\em \AddressThree}\\
        email:\hskip 4pt\texttt{\EmailThree}\\[2mm]
        \uppercase{\AuthorFour}\\
        {\em \AddressFour}\\
        email:\hskip 4pt\texttt{\EmailFour}
\or
        \uppercase{\AuthorOne}\\
        {\em \AddressOne}\\
        email:\hskip 4pt\texttt{\EmailOne}\\[2mm]
        \uppercase{\AuthorTwo}\\
        {\em \AddressTwo}\\
        email:\hskip 4pt\texttt{\EmailTwo}\\[2mm]
        \uppercase{\AuthorThree}\\
        {\em \AddressThree}\\
        email:\hskip 4pt\texttt{\EmailThree}\\[2mm]
        \uppercase{\AuthorFour}\\
        {\em \AddressFour}\\
        email:\hskip 4pt\texttt{\EmailFour}\\[2mm]
        \uppercase{\AuthorFive}\\
        {\em \AddressFive}\\
        email:\hskip 4pt\texttt{\EmailFive}
\fi
%%%%%%%%%%%%%%%%%%%%%%%%%%%%%%%%%%%%%%%%%%%%%%%%%%%%%%%%%%%%
%
%  Add Submission Acceptance Dates
%
%%%%%%%%%%%%%%%%%%%%%%%%%%%%%%%%%%%%%%%%%%%%%%%%%%%%%%%%%%%%
\vskip 6pt
\centerline{{\em Submitted\/} \Submitted, 
\/ {\em accepted in final form\/} \Accepted%
}

%%%%%%%%%%%%%%%%%%%%%%%%%%%%%%%%%%%%%%%%%%%%%%%%%%%%%%%%%%%%%%%%%%%%%%%%%%
%
%  Add Subject Classification and Keywords
%
%%%%%%%%%%%%%%%%%%%%%%%%%%%%%%%%%%%%%%%%%%%%%%%%%%%%%%%%%%%%%%%%%%%%%%%%%%

\vskip 6pt\noindent
AMS 2000 Subject classification: \SubjectClassification

%%%%%%%%%%%%%%%%%%%%%%%%%%%%%%%%%%%%%%%%%%%%%%%%%%%%%%%%%%%%%%%%%%%%%%%%%%
%
%  Abstract
%
%%%%%%%%%%%%%%%%%%%%%%%%%%%%%%%%%%%%%%%%%%%%%%%%%%%%%%%%%%%%%%%%%%%%%%%%%%

\vskip 18pt\noindent
%\emph{Abstract \vskip 2pt \Abstract }
\emph{Abstract \vskip 2pt } \Abstract

\nocite{*}
\newtheorem{theorem}{{Theorem}}[section] \newtheorem{lemma}{
  Lemma} \newtheorem{proposition}{ Proposition}[section]
\newtheorem{corollary}{ Corollary}[section]
\newtheorem{condition}{
  Condition}[section]\theoremstyle{definition}
\newtheorem{definition}{ Definition}[section]\theoremstyle{plain}
\newtheorem{defproperty}{\normalfont Property}[definition]
\newtheorem{proproperty}{}[defproperty]
\newtheorem{property}{\normalfont Property}[condition]
\newtheorem{remark}{ Remark}[section]
\renewcommand{\theequation}{\arabic{section}.\arabic{equation}}
\section{Introduction}\label{s:introduction}
In a continuous time setting, where security prices are modeled as
It\^o processes, the concept of tameness has been introduced as a
credit constrain in order to offset the so called ``doubling
strategies''.  \citename{Harrison81} \citeyear{Harrison81} and
\citename{Dybvig88} \citeyear{Dybvig88} study the role of this
constrain in ruling out doubling strategies. Generally speaking,
tameness limits the credit that an agent may have, that is used to
offset intermediate losses from trade and consumption.  This credit is
established in advance in terms of the value of money. Namely, the
credit limit is resettled every time to reflect the changes in a bank
account.  This model is a standard one in financial economics.  See
\citename{Karatzas98} \citeyear{Karatzas98}, \citename{Karatzas96}
\citeyear{Karatzas96}, and \citename{Duffie96} \citeyear{Duffie96} for
some discussion about it.  Nonetheless, in order to obtain
characterizations of non-arbitrage and completeness, strong technical
conditions are made that do not hold for very interesting models in
financial economics; see \citename{Kreps81} \citeyear{Kreps81},
\citename{Duffie86} \citeyear{Duffie86}, \citename{Back91}
\citeyear{Back91} and \citename{Hindy95} \citeyear{Hindy95} and more
recently 
\citename{Fernholz2004} \citeyear{Fernholz2004}.  Several
approaches have been taken to generalize this model.  For example,
\citename{Levental95} \citeyear{Levental95} study notions of
``arbitrage in tame portfolios'' and ``approximate arbitrage'';
\citename{Kreps81} \citeyear{Kreps81} and \citename{Delbaen94}
\citeyear{Delbaen94,Delbaen95,Delbaen95b,Delbaen96,Delbaen97,Delbaen97a,Delbaen97b,Delbaen98}
propose a notion of arbitrage called a ``free lunch''.  However these
notions are usually criticized by their lack of economic
justification.  \citename{Loewenstein00} \citeyear{Loewenstein00}
revisit the standard model of security prices as It\^o processes, and
show that the standard assumptions of positive state prices and
existence of an equivalent martingale measure exclude prices which are
viable models of competitive equilibrium and  are potentially
useful for modeling actual financial markets.  They propose the
concept of ``free snacks'' for admissible trading strategies.  Other
references are \citename{Stricker90} \citeyear{Stricker90},
\citename{Ansel92} \citeyear{Ansel92}, \citename{Delbaen92}
\citeyear{Delbaen92}, \citename{Schweizer92} \citeyear{Schweizer92},
\citename{Clark93} \citeyear{Clark93}, \citename{Schachermayer93}
\citeyear{Schachermayer93}, \citename{Lakner93} \citeyear{Lakner93} and
\citename{Willard99} \citeyear{Willard99}. 

In this paper we propose a new definition for tameness.  We call it
state tameness (see Definition $\ref{D:w-tame}$).  Loosely speaking,
we call a portfolio $\pi(t)$ a state tame portfolio if the value of
its gain process discounted by the so called ``state price density
process'' is bounded below.  For a definition of state price density
process see equation$~(\ref{E:state_price})$.  In financial terms,
this definition for tameness accounts for constrains on an agent
credit that are resettled at all times to reflect the changes in the
state of the economy.  Let us establish an analogy. In a Poker game,
it is natural to assume that the players have credit constrains,
depending on the ability of each of them to eventually cover losses.
If we think of a particular game for which one player has exhausted
his credit, but his stakes of winning are high, it is likely
that someone would be willing to take over his risk.  If the rules
of the game allowed it, this could increase his ability to obtain
credit.

We define state arbitrage, see Definition $\ref{D:Viability}$, and
characterize it.  As a consequence of Theorem \ref{T:Viability}, our
definition of non-arbitrage is an extension of non-arbitrage in the
context of standard financial markets.  See \citename{Karatzas98}
\citeyear{Karatzas98}.  Moreover, whenever  equation
$~(\ref{E:Integrability_market})$ holds and the volatility matrix is
invertible,  the existence of an equivalent martingale measure 
implies the non existence of arbitrage opportunities that are state tame,
but not conversely; see Remark \ref{R:equivalent_martingale}.  Our definition is  weaker  that the one
proposed by \citename{Levental95} \citeyear{Levental95} under the
condition that equation $~(\ref{E:Integrability_market})$ holds.  See
\citename{Levental95} \citeyear{Levental95}[Theorem 1 and Corollary
1], and \citename{Loewenstein00} \citeyear{Loewenstein00} for the
economic meaning of equation$~(\ref{E:Integrability_market})$. Our
definition of non-state arbitrage is weaker than the one proposed by
\citename{Delbaen95b} \citeyear{Delbaen95b}; see Remark~\ref{R:weaker}.  Our definition admits
the existence of ``free snacks'', see e.g., Remark
\ref{R:equivalent_martingale} and \citename{Loewenstein00}
\citeyear{Loewenstein00}[Corollary 2].  See also \citename{Loewenstein00}
\citeyear{Loewenstein00}[Corollary 2] and \citename{Loewenstein00}
\citeyear{Loewenstein00}[Example 5.3] for the economic viability of
those portfolios.

Next, we try to show the usefulness of the concept introduced.  This
is done by proving two extensions of the second fundamental theorem of
asset pricing and a theorem for valuation of contingent claims of the
American type suitable for the current context.

The question of completeness is about the ability to replicate or
access certain cash flows and not about how these cash flows are
valued.  Hence, the appropriate measure for formulating the question of
completeness is the true statistical probability measure, and not some
presumed to exist equivalent martingale probability measure.
\citename{Jarrow99} \citeyear{Jarrow99} elaborate further on this
point. We propose a valuation technique that does not require the
existence of an equivalent martingale measure and allows for pricing
contingent claims, even when the range of the volatility matrix is not
maximal.  See Theorem \ref{T:completeness}. The standard approach
relates the notion of market completeness to uniqueness of the
equivalent martingale measure; see \citename{Harrison79}
\citeyear{Harrison79}, \citename{Harrison81} \citeyear{Harrison81},
and \citename{Jarrow91} \citeyear{Jarrow91}.  \citename{Delbaen92}
\citeyear{Delbaen92} extends the second fundamental theorem for asset
prices with continuous sample paths for the case of infinitely many
assets. Other extension are \citename{Jarrow99} \citeyear{Jarrow99},
\citename{Battig99a} \citeyear{Battig99a}, and
\citename{Battig99b}\citeyear{Battig99b}. The recent paper
\citename{Fernholz2004} \citeyear{Fernholz2004} also extends valuation
theory, when an equivalent martingale measures fails to exists; they
are motivated by considerations of ``diversity''; see Remark
\ref{remark:fernholz2004} for a discussion about the connections with
this paper. 

Last, we formulate an extension of the American contingent claim
valuation theory.  See Theorem \ref{T:value_sacc}.  We provide a
valuation technique of the contingent claims of the American type in a
setting that does not require the full range of the volatility matrix.
See Theorem \ref{T:value_sacc} in conjunction with Theorem
\ref{T:completeness}.  Our approach is closer in spirit to a
computational approach.  See \citename{Karatzas88}
\citeyear{Karatzas88} and \citename{Bensoussan84}
\citeyear{Bensoussan84} to review the formal theory of valuation of
American contingent claims with unconstrained portfolios; see the
survey paper by \citename{Myneni92} \citeyear{Myneni92} as well as
\citename{Karatzas98} \citeyear{Karatzas98}. Closed form solutions are
typically not available for pricing American Options on
finite-horizons.  Although an extensive literature exist on their
numerical computation; interested readers are referred to several
survey papers and books such as \citename{Broadie96}
\citeyear{Broadie96}, \citename{Boyle96} \citeyear{Boyle96},
\citename{Carverhill90} \citeyear{Carverhill90}, \citename{Hull93}
\citeyear{Hull93}, \citename{Wilmott93} \citeyear{Wilmott93} for a
partial list of fairly recent numerical work on American Options and
comparisons of efficiency.  \setcounter{equation}{0}
\section{The model}\label{s:themodel}
In what follows we try to follow as closely as possible the notation in
\citename{Karatzas98} \citeyear{Karatzas98}, and
\citename{Karatzas96}\citeyear{Karatzas96}.  For the sake of
completeness we explicitly state all the hypotheses usually used for
financial market models with a finite set of continuous assets defined
on a Brownian filtration.  We assume a $d$-dimensional Brownian Motion
starting at $0$ $\{W(t),\mathcal{F}_t; 0 \leq t \leq T\}$ defined on a
complete probability space $(\Omega,\mathcal{F},\mathbf{P})$ where
$\{\mathcal{F}_t\}_{0\leq t\leq T}$ is the $\mathbf{P}$ augmentation
by the null sets in $\mathcal{F}^{W}_T$ of the natural filtration
$\mathcal{F}_t^W=\sigma(W(s),0\leq s\leq t)$, $0\leq t\leq T$,  and $\mathcal{F}=\mathcal{F}_T$.

We assume a risk-free rate process $r(\cdot)$, a $n$-dimensional mean
rate of return process $b(\cdot)$, a $n$-dimensional dividend rate
process $\delta(\cdot)$, a $n\times d$ matrix valued volatility
process $\left(\sigma_{i,j}(\cdot)\right)$; we also assume that $b(t)$,
$\delta(t)$, $r(t)$ and $(\sigma_{i,j}(t))$ are progressively
measurable processes.  Moreover it is assumed that
\begin{equation*}
\int_0^T (\left|r(t)\right|+ \left\|b(t)\right\|+\left\|\delta(t)\right\|+\sum_{i,j}\sigma_{ij}^2(t))\,dt < \infty
\end{equation*}
As usual we assume a bond price process $B(t)$ that evolves according
to the equation
\begin{equation}\label{E:Bond}
  dB(t)=B(t)r(t)dt,\qquad B(0)=1
\end{equation}
and $n$ stocks whose evolution of the price-per-share process $P_i(t)$
for the $i^{th}$ stock at time $t$, is given by the 
stochastic differential equation
\begin{eqnarray}\label{E:ipricestock}
  dP_i(t)=P_i(t)\left[b_i(t)dt + \sum_{1\leq j\leq d}\sigma_{ij}(t)\,dW_j(t)\right],& P_i(0)=p_i\in\left(0,\infty\right)\notag\\ &i=1,\cdots,n.
\end{eqnarray}
Let $\tau\in\mathcal{S}$ be a stopping time, where $\mathcal{S}$
denotes the set of stopping times
$\tau\colon\Omega\mapsto\left[0,T\right]$ relative to the filtration
$(\mathcal{F}_t)$.  We shall say that a stochastic process $X(t)$,
$t\in[0,\tau]$ is $(\mathcal{F}_t)$-adapted if $X(t\wedge\tau)$ is
$(\mathcal{F}_t)$-adapted, where $s\wedge t=\min\left\{s,t\right\}$,
for $s,t\in\mathbb{R}$. We consider a portfolio process
$(\pi_0(t),\pi(t))$, $t\in\left[0,\tau\right]$ to be a
$(\mathcal{F}_t)$-progressively measurable
$\mathbb{R}\times\mathbb{R}^n$ valued process, such that
\begin{eqnarray}\label{eq:local_portfolio}
\int_0^{\tau} |\sum_{0\leq i\leq n}\pi_i(t)||r(t)|\,dt\ + &
\int_0^{\tau} |\pi{\prime}(t)(b(t)+\delta(t)-r(t)\mathbf{1}_n)|\,\,dt\notag\\
& +\ \int_0^{\tau}\left\|\sigma{\prime}(t)\pi(t)\right\|^2\,dt <\infty
\end{eqnarray}
holds almost surely, with $\left\| x\right\|= (x_1^2+\cdots
+x_d^2)^{1/2}$ for $x\in\mathbb{R}^d$, and
$\mathbf{1}_n^{\prime}=(1,\cdots, 1)\in\mathbb{R}^n$. A
$(\mathcal{F}_t)$-adapted process $\left\{C(t), 0\leq t\leq
  \tau\right\}$ with increasing , right continuous paths, $C(0)=0$,
and $C(\tau)<\infty$ almost surely (a.s.)  is called a \emph{cumulative consumption
  process}.  Following the standard literature (see e.g.:
\citename{Karatzas98} \citeyear{Karatzas98},
\citename{Karatzas96}\citeyear{Karatzas96}) for a given
$x\in\mathbb{R}$ and $(\pi_0,\pi,C)$ as above, the process $X(t)\equiv
X^{x,\pi,C}(t)$, $0\leq t\leq\tau$ given by the equation
\begin{eqnarray}\label{E:wealth_process}
\lefteqn{\gamma(t)X(t)=x-\int_{\left(0,t\right]}\gamma(s)\,dC(s)}\nonumber\\
& &+\int_0^t\gamma(s)\pi^{\prime}(s)\left[\sigma(s)\,dW(s)\ + (b(s)+\delta(s)-r(s)\mathbf{1}_n))\,ds\right]
\end{eqnarray}
where $\gamma(t)$ is defined as
\begin{equation}\label{E:discount_factor}
\gamma(t)\stackrel{\varDelta}{=}\frac{1}{B(t)}=\exp\left(-\int_0^tr(s)\,ds\right),
\end{equation}
is the \emph{wealth process} associated with the initial capital $x$,
portfolio $\pi$, and cumulative consumption process $C$.
\begin{remark}\label{rq:local_portfolio}
  Let us observe that the condition defined by equation
  \eqref{eq:local_portfolio} is slightly different from the
  condition that defines a portfolio in the standard setting where the
  terminal time is not random.  In
  fact, only the former condition is needed in order to obtain a well
  defined wealth process as defined by equation
  \eqref{E:wealth_process}.
\end{remark}

We define a progressively measurable \emph{market price of risk} process 
$\theta(t)=\left(\theta_1(t),\cdots,\theta_d(t)\right)$ with values in
$\mathbb{R}^d$ for $t\in \left[0,T\right]$ as the unique process
$\theta(t)\in \ker^{\perp}(\sigma(t))$ , the orthogonal complement of
the kernel of $\sigma(t)$, such that
\begin{equation}\label{E:wviability}
  b(t)+\delta(t)-r(t)\mathbf{1}_n-proj_{\ker(\sigma^{\prime}(t))}(b(t)+\delta(t)-r(t)\mathbf{1}_n)=\sigma(t)\theta(t)\qquad\mbox{a.s.}
\end{equation}
(See \citename{Karatzas98} \citeyear{Karatzas98} for a proof that
$\theta(\cdot)$ is progressively measurable.)  Moreover, we assume that
$\theta(\cdot)$ satisfies the mild condition
\begin{equation}\label{E:Integrability_market}
  \int_0^T\left\|\theta(t)\right\|^2\,dt\:<\infty\qquad\mbox{a.s.}
\end{equation}
We define a \emph{state price density process} by
\begin{equation}\label{E:state_price}
  H_0(t)=\gamma(t)Z_0(t)
\end{equation}
where
\begin{equation}\label{E:auxiliar1}
  Z_0(t)=\exp\left\{-\int_0^t\theta^{\prime}(s)\,dW(s)\ -\frac{1}{2}\int_0^t\left\|\theta(s)\right\|^2\,ds\right\}.
\end{equation}
The name ``state price density process'' is usually given to the
process defined by equation$~(\ref{E:state_price})$ when the market is
a standard financial market; see \citename{Karatzas98}
\citeyear{Karatzas98}.  In that case the process $Z_0(t)$ is a
martingale and $Z_0(T)$ is indeed a \emph{state price density}.
However, in our setting we allow the possibility that
$\mathbf{E}Z_0(T)<1$.  \setcounter{equation}{0}
\section{State tameness and state arbitrage.  Characterization}\label{s:statetameness}
We propose the following definition for tameness.
\begin{definition}\label{D:w-tame}
  Given a stopping time $\tau\in\mathcal{S}$, a self-financed portfolio
  process $(\pi_0(t),\pi(t))$, $t\in\left[0,\tau\right]$ is said to be
  \emph{state-tame}, if the discounted gain process $H_0(t)G(t)$,
  $t\in\left[0,\tau\right]$ is bounded below, where $G(t)=G^{\pi}(t)$
  is the gain process defined as
\begin{equation}\label{E:gain_process}
G(t)=\gamma^{-1}(t)\int_0^t\gamma(s)\pi^{\prime}(s)\left[\sigma(s)\,dW(s)\ + (b(s)+\delta(s)-r(s)\mathbf{1}_n))\,ds\right].
\end{equation}
\end{definition}

\begin{definition}\label{D:Viability}
  A self finance state-tame portfolio $\pi(t)$, $t\in\left[0,T\right]$
  is said to be a \emph{state arbitrage opportunity} if
\begin{equation}\label{E:arbitrage}
  \mathbf{P}\left[H_0(T)G(T)\geq 0\right]=1,\qquad\text{and}\qquad\mathbf{P}\left[H_0(T)G(T)>0\right]>0
\end{equation}
where $G(t)$ is the gain process that corresponds to $\pi(t)$. We say
that a market $\mathcal{M}$ is state-arbitrage-free if no such
portfolios exist in it.
\end{definition}

\begin{theorem}\label{T:Viability}
  A market $\mathcal{M}$ is state-arbitrage-free if and only if the
  process $\theta(t)$ satisfies
\begin{equation}\label{E:Viability}
b(t)+\delta(t)-r(t)\mathbf{1}=\sigma(t)\theta(t)\qquad 0\leq t\leq T \mbox{ a.s.}
\end{equation}
\end{theorem}
\begin{remark}\label{R:localmartingale}
  We observe that if $\theta(t)$ satisfies
  equation$~(\ref{E:Viability})$ then for any initial capital $x$, and
  consumption process $C(t)$,

\begin{eqnarray}\label{E:localmartingale}
\lefteqn{H_0(t)X(t)+\int_{\left(0,t\right]} H_0(s)\,dC(s) }\nonumber\\
& =x+\int_0^tH_0(s)\left[\sigma^{\prime}(s)\pi_(s)-X(s)\theta(s)\right] ^{\prime}\,dW(s) &.
\end{eqnarray}
\end{remark}
\begin{proof}[Proof of Theorem \ref{T:Viability}]
  First, we prove necessity. For $0\leq t\leq T$ we define
\begin{eqnarray}
p(t)=proj_{\ker(\sigma^{\prime}(t))}(b(t)+\delta(t)-r(t)\mathbf{1}_n)\notag\\
  \pi(t)=\left\{\begin{array}{ll}
        \left\|p(t)\right\|^{-1}p(t) &\mbox{if $p(t)\neq 0$,}\\
        0 &\mbox{otherwise}
                 \end{array}\right.
               \notag
\end{eqnarray}
and define $\pi_0(t)=G(t)-\pi^{\prime}(t)\mathbf{1}_n$ where $G(t)$ is
the gain process defined by equation$~(\ref{E:wealth_process})$ with
zero initial capital, and zero cumulative consumption process. It
follows that $(\pi_0(t),\pi(t))$ is a self-financed portfolio with gain
process
\[
G(t)=\gamma^{-1}(t)\int_0^t\left\|p(s)\right\|\gamma(s)\mathbf{1}_{p(s)\neq
  0}\,ds .
\]
Since $H_0(t)G(t)\geq 0$, the non-state-arbitrage hypothesis implies
the desired result.  To prove sufficiency, assume that $\theta(t)$
satisfies equation \eqref{E:Viability}, $\pi(t)$ is a self-financed portfolio and $G(t)$ is the gain process that corresponds to
$\pi(t)$ as in Definition \ref{D:w-tame}.  Remark
\ref{R:localmartingale} implies that $H_0(t)G(t)$ is a
local-martingale.  By state-tameness it is also bounded below.
Fatou's lemma implies that $H_0(t)G(t)$ is a super-martingale.  The
result follows.
\end{proof}
\begin{remark}\label{R:non_arbitrage_in_random_time}
  We can extend the definition of state arbitrage opportunity to state
  tame portfolios defined on a random time.  It is worth to mentioning
  that Theorem \ref{T:Viability} remains true even with this
  apparently stronger definition.
\end{remark}

\begin{remark}\label{R:equivalent_martingale}
  It is well known that absence of arbitrage opportunities on tame
  portfolios is implied by the existence of an equivalent martingale
  measure under which discounted prices (by the bond price process) plus discounted cumulative
  dividends become martingales; see e.g., \citename{Duffie96}
  \citeyear{Duffie96}[Chapter 6].  If the volatility matrix
  $\sigma(\cdot)$ is invertible and equation 
  (\ref{E:Integrability_market}) holds, it is known that the non existence
  of arbitrage opportunities in tame portfolios is equivalent to
  $\mathbf{E}Z_0(T)=1.$, see e.g.,  \citename{Levental95}
  \citeyear{Levental95}[Corollary 1].  Our framework allows for the
  possibility that $\mathbf{E}Z_0(T)<1$, as is the case of, for instance, \citename{Levental95}
  \citeyear{Levental95}[Example 1].  Therefore, in the cited example,
  any arbitrage opportunity that is a tame portfolio, would not be a
  state tame portfolio.
\end{remark}
\begin{remark}\label{R:weaker}
It is known that the non existence of arbitrage opportunities in tame
portfolios implies that equation (\ref{E:Viability}) holds a.s. for
Lebesgue-almost-every $t\in\left[0,T\right]$; see e.g.  \citename{Karatzas98}
  \citeyear{Karatzas98}[Theorem 4.2].  At the same time, by
Theorem \ref{T:Viability}, non existence of arbitrage opportunities in
state-tame portfolios is equivalent to assuming that  equation
(\ref{E:Viability}) holds a.s. for Lebesgue-almost-every
$t\in\left[0,T\right]$. Under a more general setting,
\citename{Delbaen94}\citeyear{Delbaen94} have proved that
the existence of an equivalent martingale measure is equivalent to a
property called ``no free lunch with vanishing risk'' (NFLVR).  It is also
known that the concept of NFLVR is stronger that the non existence of
arbitrage opportunities in tame portfolios; see
e.g., \citename{Delbaen95b}\citeyear{Delbaen95b}[Theorem 1.3].  It
follows that our definition of non-state-arbitrage is weaker that NFLVR.
\end{remark}
\setcounter{equation}{0}
\section{State European Contingent Claims.  Valuation}\label{s:europeanvaluation}
Throughout the rest of the paper we assume that equation
~(\ref{E:Viability}) is satisfied.

A $(\mathcal{F}_t)$-progressively measurable semi-martingale
$\Gamma(t), 0\leq t\leq\tau$, where $\tau\in\mathcal{S}$ is a stopping
time is called a \emph{cumulative income process for the random time
  interval $\left(0,\tau\right]$}.  Let $X(t)$ defined by
\begin{eqnarray}\label{E:general_wealth_process}
\lefteqn{\gamma(t)X(t)=x+\int_{\left(0,t\right]}\gamma(s)\,d\Gamma(s)\ +}\nonumber\\
& &\int_0^t\gamma(s)\pi^{\prime}(s)\left[\sigma(s)\,dW(s)\ + (b(s)+\delta(s)-r(s)\mathbf{1}_n))\,ds\right],
\end{eqnarray}
where $\pi(t)$, $t\in\left[0,\tau\right]$,  is a $\mathbb{R}^n$ valued
$(\mathcal{F}_t)$-progressively measurable process 
such that
\[
\int_0^{\tau}\left(\left|\pi^{\prime}(t)(b(t)+\delta(t)-r(t)\mathbf{1}_n)\right|+\left\|\sigma^{\prime}(t)\pi(t)\right\|^2\right)\,
dt <\infty .
\]
It follows that $X(t)$ defines a wealth associated with the initial
capital $x$ and cumulative income process $\Gamma(t)$. Namely, if
$\pi_0(t)=X(t)-\pi^{\prime}(t)\mathbf{1}_n$, $(\pi_0,\pi)$ defines a
portfolio process whose wealth process is $X(t)$ and cumulative income
process is $\Gamma(t)$.  Moreover, it follows that
\begin{eqnarray}\label{E:general_localmartingale}
\lefteqn{H_0(t)X(t)-\int_{\left(0,t\right]} H_0(s)\,d\Gamma(s) }\nonumber\\
& =x+\int_{0}^tH_0(s)\left[\sigma^{\prime}(s)\pi_(s)-X(s)\theta(s)\right] ^{\prime}\,dW(s) &.
\end{eqnarray}
We say that the portfolio is \emph{state $\Gamma$-tame} if the process
$H_0(t)X(t)$ is (uniformly) bounded below.

We propose to extend the concepts of European contingent claim,
financiability and completeness.  Let $Y(t)$ $t\in\left[0,\tau\right]$
be a cumulative income process with $Y(0)=0$.  Assume that $Y$ has a
decomposition $Y(t)=Y_{loc}(t)+Y_{fv}(t)$, as a sum of a local
martingale and a process of finite variation. Let
$Y_{fv}(t)=Y_{fv}^{+}(t)-Y_{fv}^{-}(t)$ be the representation of
$Y_{fv}(t)$ as the difference of two non decreasing RCLL progressively
measurable processes with $Y^{+}_{fv}(0)=Y^{-}_{fv}(0)=0$, where
$Y^{+}_{fv}(t)$ and $Y^{-}_{fv}(t)$ are the positive and negative
variation of $Y_{fv}(t)$ in the interval $\left[0,t\right]$
respectively.  We denote by 
$\left|Y_{fv}\right|(t)=Y^{+}_{fv}+Y^{-}_{fv}(t)$ the total variation
of $Y_{fv}(t)$ on the interval $\left[0,t\right]$.  We also denote
$Y^{-}$ the process defined as $Y^{-}(t)=Y_{loc}(t)-Y^{-}_{fv}(t)$.
\begin{definition}\label{D:contingent_claim}
  Given a stopping time $\tau\in\mathcal{S}$, we shall call
  \emph{state European contingent claim (SECC) with expiration date
    $\tau$} any progressively measurable semi-martingale $Y(t)$,
  $t\in\left[0,\tau\right]$, with $Y(0)=0$, such that
  $-\int_0^{\tau}H_0(t)\,dY^{-}_{fv}(t)$ is bounded below and
\begin{equation}
\mathbf{E}\left[\int_0^{\tau}H^{2}_0(t)\, d\left\langle Y\right\rangle(t)\right]
+ \mathbf{E} \left[\int_{0}^{\tau}H_0(t)\,
  d\left|Y_{fv}\right|(t)\right] <\infty .
\end{equation}
Here $\left\langle Y\right\rangle(t)$ stands for the quadratic
variation process of the semi-martingale $Y(t)$. We define $u_{e}$ by
the formula
\begin{equation}\label{E:valuesecc}
  u_{e}=\mathbf{E}\int_{0}^{\tau}H_0(t)\, dY .
\end{equation}
\end{definition}
\begin{definition}\label{D:wcompleteness}
  A state European contingent claim $Y(t)$ with expiration date $\tau$
  is called \emph{ attainable} if there exist a state $(-Y)$-tame
  portfolio process $\pi(t)$, $t\in\left[0,\tau\right]$ with
\begin{equation}\label{E:wcompleteness}
  X^{u_e,\pi,-Y}(\tau^{-})=Y(\tau),\qquad\mbox{a.s.}
\end{equation}
The market model $\mathcal{M}$ is called \emph{state complete} if
every state European contingent claim is attainable. Otherwise it is
called \emph{state incomplete}.
\end{definition}
For the following theorem we assume $\left\{i_1<\cdots <
  i_k\right\}\subseteq \left\{1,\cdots, d\right\}$ is a set of indexes
and let $\left\{i_{k+1}<\cdots<i_{d}\right\}\subseteq
\left\{1,\cdots,d\right\}$ be its complement. Let $\sigma_i(t)$,
$1\leq i\leq k$,  be the $i^{th}$ column process for the matrix valued
process $(\sigma_{i,j}(t)), 0\leq t\leq T$.  Namely, $\sigma_i(t)$,
$1\leq i\leq k$, is the
$\mathbb{R}^n$-valued progressively measurable process whose $j^{th}$,
$1\leq j\leq d$ entry agrees with 
$\sigma_{i,j}(t)$, for $0\leq t\leq T$.  We denote by
$\sigma_{i_1,\cdots,i_k}(t)$, $0\leq t\leq T$  the $n\times k$ matrix
valued process  whose $j^{th}$ column process agrees with 
$\sigma_{i_j}(t)$, $0\leq t\leq T$ for $1\leq j\leq k$.  We shall
denote as
$\{\mathcal{F}_t^{i_1,\cdots,i_k}, 0\leq t\leq T\}$ the $\mathbf{P}$ augmentation by the
null sets of the natural filtration
$\left\{\sigma(W_{i_1}(s),\cdots,W_{i_k}(s),0\leq
s\leq t), 0\leq t\leq T\right\}$.
\begin{theorem}\label{T:completeness}
  Assume that $\theta_i(t)=0$ for $i\notin\left\{
    i_1,\cdots,i_k\right\}$, where  $\theta(t)=(\theta_1(t),
  \cdots,\theta_d(t))$ is the market price of risk.  Assume that
  $\sigma_{i_1,\cdots,i_k}(t)$ is a
  $\mathcal{F}_t^{i_1,\cdots,i_k}$-progressively measurable matrix
  valued process such that 
  $Range(\sigma_{i_{k+1},\cdots,i_{d}}(t))$
  $=Range^{\perp}(\sigma_{i_1,\cdots,i_{k}}(t))$ almost surely for
  Lebesgue-almost-every $t$. In addition assume that the interest rate
  process $\gamma$ is $\mathcal{F}^{i_1,\cdots, i_k}_t$-progressively
  measurable.  
  Then, any
  $\mathcal{F}^{i_1\cdots,i_k}_t$-progressively measurable state
  European contingent claim is attainable if and only if
  $Rank(\sigma_{i_1,\cdots,i_k}(t))=k$ a.s. for Lebesgue-almost-every $t$.  In particular, a financial
  market $\mathcal{M}$ is state complete if and only if $\sigma(t)$
  has maximal range a.s. for Lebesgue-almost-every $t$, $0\leq t\leq T$.
\end{theorem}
\begin{proof}[Proof of sufficiency]
  Let $Y(t)$, $t\in\left[0,\tau\right]$, be a
  $\mathcal{F}^{i_1,\cdots,i_k}(t)$-progressively measurable SECC with
  $\tau\in\mathcal{S}$. Define
\begin{equation}\label{E:wealth}
  X(t)=H^{-1}_0(t)\mathbf{E}\left[\int_{\left(t,\tau\right]}H_0(s)\,dY(s)\mid \mathcal{F}^{i_1,\cdots,i_k}(t)\right]\qquad\mbox{for }t\in\left[0,\tau\right] .
\end{equation}
From the representation of Brownian martingales as stochastic
integrals it follows that there exist a progressively measurable
$\mathbb{R}^d$-valued process $\varphi
^{\prime}(t)=\left(\varphi_1(t),\cdots,\varphi_d(t)\right)$,
$t\in\left[0,\tau\right]$, such that
\begin{equation}\label{E:representation}
  H_0(t)X(t)+\int_{\left(0,t\right]}H_0(s)\, dY(s)\ =u_e + \int_0^t\varphi ^{\prime}(s)\,dW(s)
\end{equation}
where $\varphi_i(t)=0$ for $i\notin\left\{i_1,\cdots,i_k\right\}$.
Define $\pi_e(t)$, $t\in\left[0,\tau\right]$, as the unique
$\mathbb{R}^n$-valued progressively measurable process such that
\begin{equation}\label{E:compatibility_portfolio}
  \sigma ^{\prime}(t)\pi_e(t)=H^{-1}_0(t)\varphi(t)+ X(t)\theta(t) .
\end{equation}
The existence and uniqueness of such a portfolio follows from the
hypotheses (see Lemma 1.4.7 in \citename{Karatzas98}
\citeyear{Karatzas98}). Define
$(\pi_e)_0(t)=X(t)-\pi(t)^{\prime}\mathbf{1}_n$.  It follows using
It\^o's formula that $X(t)$ defines a wealth process with cumulative
income process $-Y(t)$, with the desired characteristics. (To prove
the state $-Y(t)$ tameness of the portfolio $\pi_e(t)$, let  $u^{-}_e$
be the constant defined by the
equation$~(\ref{E:valuesecc})$ corresponding to the SECC $Y^{-}(t)$.  Let
$X^{-}(t)$,  $\varphi^{-}(t)$, and $\pi_e^{-}(t)$ be the  processes defined by
equations$~(\ref{E:wealth})$,$~(\ref{E:representation})$,
$~(\ref{E:compatibility_portfolio})$ respectively corresponding to the
SECC $Y^{-}(t)$; it follows that  $X(t)\geq X^{-}(t)$, $0\leq t\leq \tau$. The
$-Y(t)$ tameness of $\pi_e(t)$ is implied by the $-Y^{-}(t)$ tameness
of $\pi_e^{-}(t)$.  The latter follows by the definition of SECC.)
\end{proof}
\begin{proof}[Proof of necessity]  Let us assume that any
  $\mathcal{F}^{i_1,\cdots,i_k}_t$-progressively measurable SECC is attainable.  Let
  $f: L(\mathbb{R}^k;\mathbb{R}^n)\mapsto\mathbb{R}^k$ be a bounded
  measurable function such that: $f(\sigma)\in Kernel(\sigma)$ and 
  $f(\sigma)\neq\mathbf{0}$ if
  $Kernel(\sigma)\neq\left\{\mathbf{0}\right\}$, hold for every
  $\sigma\in L(\mathbb{R}^k;\mathbb{R}^n)$.  (See
  \citename{Karatzas96} \citeyear{Karatzas96}, p. 9).  Let us define
  $\psi(t)$ to be the bounded,
  $\mathcal{F}^{i_1,\cdots,i_k}_t$-progressively measurable process such
  that $\psi_{i_1,\cdots,i_k}=f(\sigma_{i_1,\cdots,i_k}(t))$ and
  $\psi_j(t)=0$ for $j\notin \left\{i_1,\cdots,i_k\right\}$.  We
  define the $\mathcal{F}^{i_1,\cdots,i_k}$-progressively measurable SECC by
\begin{equation}\label{E:SECC}
        Y(t)=\int_0^t\frac{1}{H_0(s)}\psi ^{\prime}(s)\,
        dW(s)\qquad\mbox{for }0\leq t\leq \tau .
\end{equation}
Let $\pi_e$ be the $-Y$ state tame portfolio with wealth process
$X^{u_e,\pi_e,-Y}$ as in equation $~(\ref{E:wcompleteness})$ and $u_e$
defined by equation$~(\ref{E:valuesecc})$.  It follows that
\begin{equation}\label{E:Martingale}
H_0(t)X^{u_e,\pi_e,-Y}(t)+\int_{\left(0,t\right]}H_0(s)\, dY(s)=u_e +\int_0^t\psi^{\prime}(s)\, dW(s)
\end{equation}
is a martingale.  Using equation~(\ref{E:general_localmartingale}),
and the representation of Brownian martingales as stochastic integrals
we obtain
\begin{eqnarray}\label{E:uniqueness}
\lefteqn{\psi_{i_1,\cdots,i_k}(t)=\sigma_{i_1,\cdots,i_k}^{\prime}(t)\pi_e(t)-X(t)\theta_{i_1,\cdots,i_k}(t)  } \nonumber\\
& \in Kernel^{\perp}(\sigma_{i_1,\cdots,i_k}(t)\cap Kernel(\sigma_{i_1,\cdots,i_k}(t))=\left\{\mathbf{0}\right\} &
\end{eqnarray}
a.s. for Lebesgue-almost-every $t$, $0\leq t\leq \tau$. The result follows.
\end{proof}
\begin{remark}\label{remark:fernholz2004}
\citename{Fernholz2004} \citeyear{Fernholz2004} are able to hedge
contingent claims of  European type when a martingale measure fails
to exists.  The framework of their paper is the same as ours, namely,
the model of security prices as It\^o processes. In addition they
assume that the eigenvalues of the stochastic $n\times n$-matrix of
variation-covariation rate processes $\sigma(t)\sigma^{\prime}(t), t\in[0,T]$
are uniformly bounded away from zero.  This latter condition implies
that equation (\ref{E:Viability}) holds;  as a consequence their results on valuation are implied by Theorem \ref{T:completeness}. 
\end{remark}
\setcounter{equation}{0}
\section{State American Contingent Claims.  Valuation.}\label{s:valuationamerican}
\begin{definition}\label{D:SACC}
  Let $(\Gamma(t),L(t))$, $0\leq t\leq \tau$, a couple of RCLL
  progressively measurable semi-martingales where $\Gamma(t)$,
  $t\in\left[0,\tau\right]$, is a cumulative income process with
  $\Gamma(0)=0$.  Assume that the process
\begin{equation}\label{E:discounted_payoff_process}
Y(t)=\int_{\left(0,t\right]}H_0(s)\,d\Gamma(s)+L(t)H_0(t)\qquad\mbox{for
}0\leq t\leq\tau ,
\end{equation}
is a continuous semi-martingale such that $Y$ and $L(t)H_0(t)$, $0\leq
t\leq \tau$, are uniformly bounded below.
We shall call a \emph{state American contingent claim (SACC)} a couple
of processes as above such that
\begin{equation}\label{E:value_sacc}
u_a=\sup_{\tau^{\prime}\in\mathcal{S}(\tau)}\mathbf{E}[Y(\tau^{\prime})]<\infty
,
\end{equation}
where
$\mathcal{S}(\tau)=\left\{\tau^{\prime}\in\mathcal{S};\tau^{\prime}\leq\tau\right\}$.
We shall call the process $Y(t)$ the \emph{discounted payoff process},
$L(t)$ the \emph{lump-sum settlement process} and $u_a$ the
\emph{value of the state American contingent claim}.
\end{definition}
\begin{theorem}\label{T:value_sacc}
  Let $\left\{i_1,\cdots,i_k\right\}\subseteq
  \left\{1,\cdots,d\right\}$ be a set of indexes. Assume the
  hypotheses of theorem~\ref{T:completeness}.  If $(\Gamma(t),L(t))$
  is a state American contingent claim where the discounted payoff
  process is $\mathcal{F}^{i_1,\cdots,i_k}_t$-progressively measurable then there
  exist a $-\Gamma(t)$ state tame portfolio $\pi_a$ such that
\begin{equation}\label{E:feasability_sacc}
X^{u_a,\pi_a,-\Gamma}(t) \geq L(t)\qquad\mbox{a.s. for  }0\leq t\leq
\tau .
\end{equation}
Indeed,
\begin{eqnarray}\label{E:optimality_sacc}
u_a=\inf\{u\in\mathbb{R}\mid&\text{there exist a }-\Gamma(t) \text{ state tame portfolio }\notag\\
& \pi\text{ with }X^{u,\pi,-\Gamma}(t)\geq L(t)\text{ a.s. for }0\leq
t\leq \tau \} .
\end{eqnarray}
\end{theorem}
\begin{lemma}\label{lem:american_aproximation}
  Given $\tau_1, \tau_2\in\mathcal{S}\left(\tau\right)$, there exist
  $\tau^{\prime}\in\mathcal{S}\left(\tau\right)$ with
\[
u_a\geq\mathbf{E}\left[Y\left(\tau^{\prime}\right)\right]\geq\max\left\{\mathbf{E}\left[Y\left(\tau_1\right)\right],\mathbf{E}\left[Y\left(\tau_2\right)\right]\right\}
\]
  such that
\[
\mathbf{E}\left[Y\left(\tau^{\prime}\right)\mid
  \mathcal{F}_t\right]\geq\max\left\{ \mathbf{E}\left[Y(\tau_1)\mid
    \mathcal{F}_t\right],
  \mathbf{E}\left[Y\left(\tau_2\right)\mid\mathcal{F}_t\right]\right\}\mbox{
  for all $t\in \left[0,\tau\right]$}. 
\]
\begin{proof}
  Define
\begin{eqnarray*}
    \label{eq:1}
\lefteqn{\tau^{\prime}=\tau_1\wedge\tau_2\mathbf{1}_{\mathbf{E}\left[Y\left(\tau_1\vee\tau_2\right)\mid\mathcal{F}_t\right]\left(\tau_1\wedge\tau_2\right)<Y\left(\tau_1\wedge\tau_2\right)}}\\
& &\mbox{} + \tau_1\vee\tau_2\mathbf{1}_{\mathbf{E}\left[Y\left(\tau_1\vee\tau_2\right)\mid\mathcal{F}_t\right]\left(\tau_1\wedge\tau_2\right)\geq
  Y\left(\tau_1\wedge\tau_2\right)} 
  \end{eqnarray*}
  where $s\vee t=\max\{s,t\}$, and $s\wedge t=\min\{s,t\}$.  Then
  $\tau^{\prime}$ has the required properties.
\end{proof}
\end{lemma}
\begin{proof}[Proof of Theorem \ref{T:value_sacc} ]
  Let $Y(t)$, $0\leq t\leq\tau$, be the discounted payoff process.  There
  exist a sequence of stopping times $\left(\sigma_n\right)$ in
  $\mathcal{S}\left(\tau\right)$ such that
  $\mathbf{E}\left[Y\left(\sigma_n\right)\right] \uparrow u_a$,
  $\mathbf{E}\left[Y\left(\sigma_{n+1}\right)\mid\mathcal{F}_{t}\right]\geq\mathbf{E}\left[Y\left(\sigma_n\right)\mid\mathcal{F}_t\right]$
  for $t\in\left[0,\tau\right]$, with the property that for any
  rational $q\in\mathbb{Q}\cap\left[0,T\right]$, there exist
  $N_q\in\mathbb{N}$ such that
  $\mathbf{E}\left[Y\left(\sigma_n\right)\mid\mathcal{F}_t\right]\left(q\wedge\tau\right)\geq
  Y\left(q\wedge\tau\right)$ . The latter follows by lemma
  \ref{lem:american_aproximation}.  By Doob's inequality, $\mathbf{E}\left[Y(\sigma_n)\mid\mathcal{F}_t\right]$
  is a Cauchy sequence in the sense of uniform convergence in
  probability.  By completeness of the space of local-martingales,
  there exist a local-martingale $\overline{Y}(t)$,
  $t\in\left[0,\tau\right]$, such that
  $\mathbf{E}\left[Y(\sigma_n)\mid\mathcal{F}_t\right]\rightarrow\overline{Y}(t)$, $t\in\left[0,\tau\right]$, uniformly in probability.  It
  follows by continuity that $\overline{Y}(t)\geq Y(t)$ for
  $t\in\left[0,\tau\right]$, and clearly $\overline{Y}(0)=u_a$.
  Define $\tau_n$ to be the first hitting time of $\overline{Y}(t)$, $t\in\left[0,\tau\right]$, to the set
  $\left[-n,n\right]^c$.  From the representation of Brownian
  martingales as stochastic integrals it follows that there exist a
  progressively measurable $\mathbb{R}^d$-valued process $\varphi
  ^{\prime}(t)=\left(\varphi_1(t),\cdots,\varphi_d(t)\right)$,
  $t\in\left[0,\tau_n\right]$, such that
\begin{equation}\label{E:local_representation}
  \overline{Y}(t)=u_a + \int_0^t\varphi ^{\prime}(s)\,dW(s)
\end{equation}
where $\varphi_i(t)=0$ for $i\notin\left\{i_1,\cdots,i_k\right\}$.
Define $X(t)$, $t\in\left[0,\tau\right]$, by
\begin{equation*}
   H_0(t)X(t)+\int_{\left(0,t\right]}H_0(s)\, d\Gamma(s)\
   =\overline{Y}(t) .
\end{equation*}
Define $\pi_a(t)$, $t\in\left[0,\tau\right]$, as the unique
$\mathbb{R}^n$-valued progressively measurable process such that
\begin{equation*}
  \sigma ^{\prime}(t)\pi_a(t)=H^{-1}_0(t)\varphi(t)+ X(t)\theta(t) .
\end{equation*}
The existence and uniqueness of such a portfolio follows by the
hypotheses (see Lemma 1.4.7 in \citename{Karatzas98}
\citeyear{Karatzas98}). Define
$(\pi_a)_0(t)=X(t)-\pi_a(t)^{\prime}\mathbf{1}_n$.  It follows using
It\^o's formula that $X(t)$ defines a wealth process with cumulative
income process $-\Gamma(t)$, $t\in\left[0,\tau\right]$, with the
desired characteristics.  Equation~\eqref{E:optimality_sacc} is a
consequence to the fact that the discounted payoff process is a
super-martingale.
\end{proof}
\begin{remark}\label{rq:new_portfolio}
  Let us observe that it is not possible to obtain optimal stopping times for the
  version of the theorem for valuation of American contingent claims
  that we presented.  Nonetheless, it is worth to point out that the
  conditions of the Theorem \ref{T:value_sacc}, are probably the
  weakest possible.
\end{remark}

\section[Acknowledgements]{Acknowledgements}
I thank Professor  J. Cvitani\'c, and Professor N. E. Gretsky
for suggestions made on a preliminary draft of this paper.  I also want
to thank professor  M. M. Rao for a detailed reading of the first
version of this paper and suggestions made on it that led to a
substantial improvement of the paper,  an anonymous referee for
valuable suggestions, and an associate editor for pointing out the
recent paper \citename{Fernholz2004} \citeyear{Fernholz2004} and its
connections with this work.
\ifx\undefined\BySame
\newcommand{\BySame}{\leavevmode\rule[.5ex]{3em}{.5pt}\ }
\fi
\ifx\undefined\textsc
\newcommand{\textsc}[1]{{\sc #1}}
\newcommand{\emph}[1]{{\em #1\/}}
\let\tmpsmall\small
\renewcommand{\small}{\tmpsmall\sc}
\fi

\end{document}